\documentclass[12pt, reqno]{amsart}

\usepackage{amssymb,latexsym,amsmath,amsfonts}
\usepackage{mathrsfs}
\usepackage{graphicx}
\usepackage[usenames]{color}
\usepackage{hyperref}
\usepackage{comment}
\usepackage{enumitem}

\definecolor{DPurple}{rgb}{0.46,0.2,0.69}

\numberwithin{equation}{section}

\allowdisplaybreaks

\theoremstyle{definition}

\theoremstyle{remark}

 \theoremstyle{plain}



\begin{document}

\title[On some new metric characterisations of Hilbert spaces]{On some new metric characterisations of Hilbert spaces}

\parindent=0mm \vspace{.2in}
\author{M.A. Sofi}
\address{Department of Mathematics, University of Kashmir, Srinagar-190006, India.}
\email{aminsofi@gmail.com}


\begin{abstract}
In the literature surrounding the theory of Banach spaces, considerable effort has been invested in exploring the conditions on a Banach space X that characterise X as being an inner product space or as a linearly isomorphic copy of a Hilbert space. On the other hand, a different theory emerges when the class of Banach spaces is looked upon as a Lipschitz category where Lipschitz maps are used as morphisms in the new category in place of the familiar bounded linear maps in the linear theory.

\parindent=8mm \vspace{.1in}
This paper provides a short survey of recent results involving the appropriate Lipschitz analogues of certain well known results from the linear theory characterizing Hilbert spaces. Whereas isometric description of Hilbert spaces has all along been a popular theme in this line of investigations, we shall concentrate mainly on isomorphic characterisations which entail the existence of an equivalent norm on the underlying space arising from an inner product. 

\end{abstract}
\maketitle

\section{Introduction and motivation}\label{S:intro}

\parindent=0mm \vspace{.0in}
Ever since the inception of Banach space theory about a century ago, it has been a favourite theme in functional analysis to identify Banach space properties that characterise the norm of the space arising from an inner product. The study of such properties belongs in the isometric theory of Banach spaces whereas the isomorphic theory involves the investigation of conditions that characterise when the norm of the space is equivalent to another norm that is induced by an inner product. In other words, the given Banach is linearly isomorphic with a Hilbert space.

\parindent=0mm \vspace{.1in}
On the other hand, it turns out that the (nonlinear) metric induced by the norm on a Banach space is so tightly tied to the linear structure of a Banach space that if two Banach spaces are in bijective correspondence via a (nonlinear) isometry or a biLipschitz map, then their linear structures are preserved to a very great degree. An example of this phenomenon is provided by the famous Mazur-Ulam theorem [16] which says that as soon as two (real) normed spaces are isometric via an isometry $T$ as metric spaces, then $T$ is already linear if it fixes the origin. On the other hand, a celebrated theorem of P. Enflo [5] guarantees that if a Banach space is uniformly homeomrphic with a Hilbert space, then it is linearly isomorphic with it.

\parindent=0mm \vspace{.1in}
The above discussion reveals that by virtue of the metric induced by the norm on a Banach space, there are weaker structures attached to it which arise from the uniform structure or the coarse structure of the induced metric, with the associated morphisms consisting of uniformly continuous or coarse continuous maps. On the other hand, it turns out that apart from the (metric) isometries, the class of Lipschitz maps provides the appropriate morphisms in the new setting of the Lipschitz category of Banach spaces.

\parindent=0mm \vspace{.1in}
The present work is devoted to the description (without proofs, except in some important cases where main ideas of the proof are sketched) of a brief summary of the results involving Lipschitz characterization of Hilbert spaces which have been proved in the last decade including some results of the author which have already appeared in the literature apart from those where work is still in progress. To this end, we shall be mainly concerned with the Lipschitz analogues of the following results which are well known from the linear structure theory of Banach spaces. More precisely, we shall explore the extent to which the Lipschitz analogues of the following phenomena hold true.

\parindent=0mm \vspace{.1in}

(a).	Lindenstrauss Tzafiriri Complemented Subspaces Theorem.

\parindent=0mm \vspace{.1in}
(b)(i). Grothendieck’s Theorem on linear embedding of a Banach space and its dual into $L_1$-spaces.

\parindent=0mm \vspace{.05in}
(ii). Linear embedding of a Banach space into $L_p$-spaces.

\parindent=0mm \vspace{.05in}
(c).  Existence of bounded linear extension operators on the space of bounded linear functions on subspaces of a Banach space.

\parindent=0mm \vspace{.05in}
(d). Linear selection involving norm attaining functionals on a Banach space.

\parindent=0mm \vspace{.1in}
Before we describe these issues, we list below some examples of situations involving nonlinear phenomena which characterise Hilbert spaces (linearly isometrically).

\parindent=0mm \vspace{.1in}
{\bf Example.1.1 [2]:} Given a real Banach space $X$, then each of the following conditions is equivalent to $X$ being Hilbertian.

\parindent=0mm \vspace{.1in}
(a).~~$\left\| f\right\|\,\left\| g\right\|\le 2\left\| fg\right\|,\quad \forall~f,g\in X^{*}$.

\parindent=0mm \vspace{.1in}
(b).~~$\|L\|=\|L^{\wedge}\|$, for each symmetric bilinear form $L$ on $X$. (Here $L^{\wedge}$ denotes the polynomial associated to $L$). 

\parindent=0mm \vspace{.1in}
(c).~~$\forall ~x(\neq 0),y(\neq 0)\in X,$ there exist a symmetric bilinear form $L$ on $X$ such that $L(x,y)=\|x\|\,\|y\|$.

\parindent=0mm \vspace{.1in}
{\it Comments} (i). The norm of the product in the RHS of $(a)$ is meant in the sense of $fg$ acting as a bilinear form on $X: fg(x,y)=f(x)g(y),x,y\in X$. In general, the norm of an $n$-linear $L$ on $X$ is defined by: 
\begin{align*}
\big\|L\big\|=\sup \left\{\big|L(x_1,x_2,\dots,x_n)\big|:x_1,x_2,\dots,x_n\in B_X \right\}.
\end{align*}
Also $L^{\wedge}$ in $(b)$ is defined to be the (polynomial) mapping: $L^{\wedge}(x)=L(x,x,\dots,x).$

\parindent=0mm \vspace{.1in}
(ii). A complex analogue of $(c)$ is also valid whereby $L$ is required to be Hermitian:  $L(x,y)=\overline{L(y,x)},\forall~ x,y\in X.$

\parindent=0mm \vspace{.1in}
(iii). Without requiring symmetry of $L$, the assertion in Ex.1(c) always holds in an arbitrary normed space by choosing $L(x, y)= f(x)g(y)$. This follows from the Hahn Banach extension theorem. 

\parindent=0mm \vspace{.1in}
{\bf Example.1.2 (Foias, see [20], p.27):} Given a (complex) Banach space $X$ such that for each contraction $T$ on $X$ and a polynomial $p(z)$, it holds that
\begin{align*}
\big\|p(T)\|\le \sup_{\|z\|\le 1} \big|p(z)\big|,
\end{align*}
then $X$ is linearly isometric to a Hilbert space.

\parindent=0mm \vspace{.1in}
{\bf Example.1.3 ([13]):} Recall that given a closed subspace $M$ of a Hilbert space $H$, then we can write $H=M\bigoplus M^{\perp}$ where $M^{\perp}$, the orthogonal complement of $M$ is a closed subspace of $H$. A celebrated theorem of Lindenstrauss and Tzafriri [13] shows that this is a characteristic property of Hilbert spaces. Precisely, if $M$ is a closed subspace of a Banach space $X$ such that there exists a closed subspace $N$ of $X$ with $X=M\bigoplus N$, then $X$ is a linearly isomorphic to a Hilbert space. As we shall see, this theorem and some of its consequences will motivate most of what is to follow in this paper (See Section 3.A).

\section{Banach spaces as a Lipschitz category }

{\bf Definition 2.1:} A map $f:M\to N$ acting between metric spaces $M$ and $N$ is said to be a Lipschitz map if for some $c\ge 1$, the following holds
\begin{align*}
\big\|f(x)-f(y)\big\|\le c\big\|x-y\big\|,\quad x,y\in M.
\end{align*}

\parindent=0mm \vspace{.0in}

For $c=1$, we say that $f$ is a non-expansive mapping. The space of all Lipschitz maps $f:M\to R$ shall be denoted by $Lip(M)$, whereas we shall mostly consider the (sub) space $Lip_0(M)$ of $Lip(M)$ consisting of functions vanishing at a distinguished point, say $\theta\in M$. It is easily checked that under pointwise operations, $Lip_0(M)$ is a Banach space (the Lipschitz dual of $M$) when equipped with the norm: 
\begin{align*}
\big\|f\big\|=\sup_{x\neq y} \dfrac{\big|f(x)-f(y)\big|}{d(x,y)}.
\end{align*}

\parindent=0mm \vspace{.0in}
As already indicated, the term Lipschitz category shall be used for the category comprised of Banach spaces as its objects and Lipschitz maps acting between them as morphisms.

\parindent=0mm \vspace{.1in}
We state below the Lipschitz counterpart of the well-known Hahn Banach extension from linear functional analysis.

\parindent=0mm \vspace{.1in}
{\bf Theorem 2.2} (E. J. McShane, [18]): Every Lipschitz function on a subset of a metric space $M$ can be extended to a Lipschtiz function on $M$.

\parindent=0mm \vspace{.1in}
For a comprehensive treatment of Lipschitz functions and the theory surrounding the spaces of Lipschitz functions and their isometric preduals, the so-called Lipschitz free spaces, we recommend the encyclopaedic treatise [3].

\section{Linear versus Lipschitz theory}

\parindent=0mm \vspace{.1in}
{\bf (A). A glimpse into the linear theory of Banach spaces and its Lipschitz counterpart}

\parindent=0mm \vspace{.1in}
To motivate what is to follow next, we begin by listing below four important theorems on the characterisation of Hilbert space which are well known from the linear structure theory of Banach spaces.

\parindent=0mm \vspace{.1in}
{\it Complemented subspaces Theorem}

\parindent=0mm \vspace{.1in}

{\bf Theorem 3.A.1} (Lindenstrauss and Tzafiriri [13]); A Banach space $X$ is linearly isomorphic with a Hilbert space if and only if all of its closed subspaces are complemented in $X$. 

\parindent=0mm \vspace{.1in}
Some important consequences of this important statement (LT) are provided below.

\parindent=0mm \vspace{.1in}
{\bf Theorem 3.A.2:} For a Banach space $X$ with $\dim X>2$, TFAE:

\parindent=0mm \vspace{.1in}
(i). A bounded linear map defined on a subspace of X and taking values in an arbitrary Banach space Z extends to a bounded linear map on $X$.

\parindent=0mm \vspace{.1in}
(ii). Every closed subspace of $X$ is complemented in X. 

\parindent=0mm \vspace{.1in}
(iii). $X$ is (isomorphically) a Hilbert space.

\parindent=0mm \vspace{.1in}
An isometric counterpart of the previous theorem based on an old theorem of Kakutani [9] and independent of the complemented subspace theorem is given by the following theorem. 

\parindent=0mm \vspace{.1in}

{\bf Theorem 3.A.3:} For a Banach space $X$ with $\dim X>2,$ TFAE:

\parindent=0mm \vspace{.1in}
(i). A bounded linear map defined on a subspace of $X$ and taking values in an arbitrary Banach space $Z$ extends to a bounded linear map on $X$ without increase of norm.

\parindent=0mm \vspace{.1in}
(ii) Every closed subspace of $X$ is 1-complemented in $X$. 

\parindent=0mm \vspace{.1in}
(iii). $X$ is (isometrically) a Hilbert space.

\parindent=0mm \vspace{.1in}
Another important consequence of (LT) is the following beautiful characterisation of Hilbert space in terms of the existence of a bounded linear selection of Hahn Banach extensions.

\parindent=0mm \vspace{.1in}
{\bf Theorem 3.A.4([10], [25]):} For a Banach space $X$, TFAE:

\parindent=0mm \vspace{.1in}
(i). For each closed subspace $Y$ of $X$, there exists a bounded linear “extension map” $G: Y^{*}\to X^{*}: G(g)(y) = g(y), g\in Y^{*},y\in Y$.

\parindent=0mm \vspace{.1in}
(ii). $Y^{**}$ is complemented in $X^{**}$.

\parindent=0mm \vspace{.1in}
(iii) $X$ is a Hilbert space.

\parindent=0mm \vspace{.1in}
{\bf Remark 3.A.5:} It turns out that the properties as indicated in the previous theorem are equivalent to $Y$ being locally complemented in $X$ in the following sense:

\parindent=0mm \vspace{.1in}
{\bf Definition 3.A.6:} A subspace $Y$ of a Banach space $X$ is said to be locally $c$-complemented if there exists $c>0$ such that for each finite dimensional subspace $M$ of $X$, there exists a continuous linear map $f:M\to Y$ with $\|f\|\le c$ and $f(x)=x$ for all $x\in M\cap Y$.

\parindent=0mm \vspace{.1in}
In a work [25] that is in progress, the author has proposed a different approach to a proof of Theorem A.1.4 which makes use of an equality involving numerical parameters connecting the projection constant of $X$ and the minimum norm of the extension operator arising from Theorem A.1.4. Thus, given a subspace $Y\subset X$, we denote
\begin{align*}
E_{Y}(X)&=\inf \left\{\|T\|: T: Y^{*}\to X^{*}~\text{  is an extension operator} \right\}\\
P_{Y}(X)&=\inf \left\{\|P\|: P: X\to Y~\text{  is a projection} \right\}\\
E(X)&=\sup\left\{E_Y(X): \dim Y<\infty \right\}\\
P(X)&=\sup\left\{P_Y(X): \dim Y<\infty \right\}.
\end{align*}

\parindent=0mm \vspace{.0in}

We have

\parindent=0mm \vspace{.1in}
{\bf Theorem 3.A.7 ([25]):} Given a Banach space $X$, and assume that each (closed) subspace $Y$ of $X$ is locally complemented. Then $P(X)=E(X)<\infty$.

\parindent=0mm \vspace{.1in}
The following corollary provides a strengthening of the Lindenstrauss-Tzafiriri  complemented subspace theorem. 

\parindent=0mm \vspace{.1in}
{\bf Corollary 3.A.8:} Every closed subspace of a Banach space $X$ is locally complemented if and only if it is a Hilbert space.

\parindent=0mm \vspace{.1in}
{\bf Proof:} By the previous theorem, $P(X)=E(X)=c<\infty.$ Using a well-known technique involving finite dimensional spaces, this gives: $\sup_{E\subset X} d(E,\ell_2^{\dim E})\le 4 c^2$. Finally, by Joichi’s theorem [8], it follows that $X$ is (isomorphic to)  a Hilbert space.

\parindent=0mm \vspace{.1in}
{\bf Remark 3.A.9:} The above theorems serve to unify all the previously known results involving the extendibility of maps surrounding the characterization of  Hilbert and $L_1$-predual spaces and the fundamental duality between these classes of spaces.  It turns out that the basic idea underlying this duality is rooted in the role of local complementedness in describing the two classes of spaces. Combining Theorems 2.2 and 2.3, it follows that a Banach space $X$ is

\parindent=0mm \vspace{.1in}
1. An $L_1$-predual space if and only if it is locally complemented in each Banach space containing $X$.

\parindent=0mm \vspace{.1in}
2. A Hilbert space (isomorphically) if  and only if each of its closed subspaces is locally complemented in $X$.

\parindent=0mm \vspace{.1in}
Yet another interesting consequence of the (CST) is the following characterization of Hilbert space in terms of the existence of a nice right inverse of the quotient map. 

\parindent=0mm \vspace{.1in}
{\bf Theorem 3.A 10:} Given a Banach space $X$ such that for each closed subspace $M$ of $X$ and the quotient map $\varphi X\to X/M$, there exists a bounded linear right inverse $\psi:X/M\to X~(\varphi\psi = \text{identity map on}\, X/M)$, then $X$ is (linearly and isomorphically) a Hilbert space.

\parindent=0mm \vspace{.1in}
The proof follows by noting that the existence of $\psi$ as indicated above is equivalent to M being complemented in $X$.

\parindent=0mm \vspace{.1in}

{\bf (B). In search of a Lipschitz analogue of the Complemented Subspaces Theorem}

\parindent=0mm \vspace{.1in}
To this end, we follow an approach which is different from one that is used in the linear theory. Our main tool will be the existence of bounded extension operators between spaces of Lipschitz functions defined on Banach spaces and their subspaces and derive from that an appropriate Lipschitz analogue of the (CST) and its consequences as described above.

\parindent=0mm \vspace{.1in}
Bounded Extension operators on spaces of Lipschitz functions

\parindent=0mm \vspace{.1in}
{\bf Definition 3.B.1(i):} Given metric spaces $M, N$ and function spaces $F(A, N), A\subseteq M$, we say that a map $\psi: F(A,N)\to F(M,N)$ is an extension operator  if $:\psi(g)\big|_A=g$ for each $g\in F(A,N)$.

\parindent=0mm \vspace{.1in}
(i). A metric space is said to satisfy the ``Lipschitz Selection Property" (LSP) if the map $\psi :Lip(A)\to Lip(M):\psi(g)\big|_A=g$, for each $g\in Lip(A)$ is an extension operator for all subsets $A$ of $M$. 

\parindent=0mm \vspace{.1in}
In other words, (LSP) entails the existence of an linear ``extension" operator $\psi: Lip(A)\to Lip(M)$  such that $\psi(g)\big|_A=g$, for each $g\in Lip(A)$ where $A$ is chosen to be an arbitrary subset of the metric space $M$.

\parindent=0mm \vspace{.1in}
{\bf Problem:} Whether the choice of a Lipschitz extension can be made linearly and continuously.

\parindent=0mm \vspace{.1in}
{\bf Example 3.B.2:} The following metric spaces have the (LSP):

\parindent=0mm \vspace{.1in}
(i). $\mathbb R^n$ (with respect to any norm).

\parindent=0mm \vspace{.1in}
(ii). Heisenberg group.

\parindent=0mm \vspace{.1in}
Considering the Lipschitz dual of a Banach space in place of its linear dual, the following theorems provide appropriate analogue of Theorem A.1.4 in the Lipschitz category.

\parindent=0mm \vspace{.1in}
{\bf Theorem 3.B.3([21]):} Let $X$ be a Banach space and $Z$ a subspace of $X$ such that there exists a bounded linear map ``extension" map $F:Lip_0(Z)\to Lip_0(X)$. Then $X$ is Hilbert, and conversely.

\parindent=0mm \vspace{.1in}
The proof follows from Corollary 3.A.8 combined with the following:

\parindent=0mm \vspace{.1in}
{\bf Theorem 3.B.4 ([24]):} Let $X$ be a Banach space and $Z$ a subspace of $X$ such that there exists a bounded linear map ``extension" map $F:Lip_0(Z)\to Lip_0(X)$. Then there exists a bounded linear extension map $G: Z^{*}\to X^{*}$. Hence $Z$ is c-locally complemed with $\|F\|\ge c$.

\parindent=0mm \vspace{.1in}
The last assertion is a well known theorem of Kalton [10].

\parindent=0mm \vspace{.1in}
{\bf Corollary 3.B.5:} For a Banach space $X$, TFAE:

\parindent=0mm \vspace{.1in}
(i). $X$ is a Hilbert space.

\parindent=0mm \vspace{.1in}
(ii). Each closed convex subset of $X$ is a Lipschitz retract.
	
\parindent=0mm \vspace{.1in}	
(iii). For each (closed) subspace $Z$ of $X$, there exists a bounded linear ``extension" map $F:Lip_0 (Z)\to Lip_0(X)$.

\parindent=0mm \vspace{.1in}
(iv). For each (closed) subspace $Z$ of $X$, there exists a bounded linear ``extension" map  $G: Z^*\to X^*$.

\parindent=0mm \vspace{.1in}
(v). Each (closed) subspace of $X$ is locally complemented.

\parindent=0mm \vspace{.1in}
Back to Theorem 3.B.3 where we have seen how the existence of a bounded linear extension map $F:Lip_0(Z)\to Lip_0(X)$ for each subspace $Z$ of $X$ yields $X$ as a Hilbert space. In this context, a natural question arises whether the conclusion of Theorem 3.B.3 still holds if $F$ is assumed to be a Lipschitz extension map.  

\parindent=0mm \vspace{.1in}
We show that for subspaces $Z$ of $X$ for which $Lip_0(Z)$ is separable, the answer is in the affirmative. In particular, this holds if $X$ is assumed to be separable. This follows from an important theorem of Godefroy and Kalton [6] to the effect that given Banach spaces $X, Z$, a separable and a bounded linear operator $T:X\to Z$, the existence of a Lipschitz right inverse $S:Z\to X$ of T yields the existence of a bounded linear right inverse $L:Z\to X$ of $T:To L=Id_Z$.

\parindent=0mm \vspace{.1in}
Indeed, consider the bounded linear map $\psi: Lip_0(X)\to Lip_0(Z)$ given by: $\psi f=f/Z$. Now the presumed existence of $F:Lip_0(Z)\to Lip_0(X)$  as a Lipschitz extension operator gives that $F$ is a Lipschitz right inverse of $\Psi$. By the Godefroy-Kalton theorem just stated, we get a bounded linear right inverse $G$ of  $\Psi$ i.e., $G:Lip_0(Z)\to Lip_0(X)$ is a bounded linear map such that $\Psi o G=Id_{Lip_0(Z)}$. Finally, $G$ is an extension map: given $g\in Lip_0(Z)$, we have $\Psi o G=(g)=g$ and so, $G(g)\big|_Z=g$. In other words, $G$ is a bounded linear extension operator and so we are in the situation of Theorem 3.B.3 to arrive at the desired conclusion. 

\parindent=0mm \vspace{.1in}
We now provide a sample of recent results involving the implications of the existence of nice liftings of a quotient map on a large collection of quotient spaces.

\parindent=0mm \vspace{.1in}
The situation where $\varphi:X\to X/M$ has a Lipschitz inverse for each closed subspace $M$ of $X$ is covered by the following theorem:

\parindent=0mm \vspace{.1in}
{\bf Theorem 3.B.6([24]):} Given a Banach space $X$ such that for each closed
subspace $M$ of $X$ and the quotient map $\varphi:X\to X/M$, there exists a Lipschitz map $\psi:X/M\to X$  such that  $\varphi\psi= \text{identity map on}\, X/M$, then $X$ is (isomorphically) a Hilbert space.

\parindent=0mm \vspace{.1in}
Under the conditions of the theorem, it follows that $M$ is a Lipschitz retract under the map: $r(x)=x-\psi \varphi(x)+\psi(0)$ and apply apply Corollary 3.B.5  $((ii)\Longrightarrow(i))$.

\parindent=0mm \vspace{.1in}
{\bf Corollary 3.B.7:} If $Z$ is a subspace of $X$ which is a $c$-Lipschitz retract and is complemented in its bidual, then $Z$ is $c$-complemented in $X$.

\parindent=0mm \vspace{.1in}
{\bf Proof:} Let $P: Z^{**}\to Z$ be a bounded linear projection and let $G:Lip_0(Z)\to Lip_)(X)$ be an ``extension" operator induced by a Lipschitz retraction of $X$ onto $Z$. As shown in the proof of Theorem 3.B.4, (the restriction of) $G$ may be considered as an extension operator $G: Z^{*}\to X^{*}$. Note that the adjoint map $G^{*}: X^{**}\to Z^{**}$ is such that for $z\in Z\subseteq X\subseteq X^{**}, G^{*} (z^{\circ})=z$ where $z$ is identified with the evaluation map $z^{\circ}\in X^{**}$. Indeed, for  since $Q(x)\in Z$, it follows that $Q^2(x)=\left(P\left(G^*\big|_X\right)\right)Q(x) =P(Q(x))=Q(x)$.

\parindent=0mm \vspace{.1in}
A local theoretic analogue of Theorem 3.B.5 also holds by combining Corollary 3.A.8 with  the famous Lindenstrauss-Tzafiriri theorem which asserts that a Banach space $X$ is (isomorphically) a Hilbert space if (and only if) its finite dimensional subspaces are uniformly complemented.

\parindent=0mm \vspace{.1in}
{\bf Theorem 3.B.8 ([25]):} Assume that for each finite dimensional subspace $M$ of $X$, there exists $c>0$ such that the quotient map $\varphi: X \to X/M$, admits a Lipschitz right inverse $\psi: X/M\to X:\varphi\psi=\text{identity map on }\, X/M$  with $\|\psi\|\le c$. Then $X$ is (isomorphically) a Hilbert space.

\parindent=0mm \vspace{.1in}
Indeed, as note in Theorem 3.B.6, the given condition yields that each finite dimensional subspace of $X$ is a $(c+1)$-Lipschitz retract, and hence $(c+1)$-locally complemented. Combining this observation with Theorem 3.A.7 gives that finite dimensional subspaces of $X$ are $(c+1)$-complemented and the (proof of the) Lindenstrauss-Tzafiriri theorem completes the argument.

\parindent=0mm \vspace{.1in}
Using an idea of Kalton [11], it's possible to provide the following stronger
version of Theorem 3.B.6. We’ll need the following definition.

\parindent=0mm \vspace{.1in}
{\bf Definition 3.B.9:} A subset $G$ of a metric space $(M, d)$ is called a proper if there exist $a,b>0$ such that $d(x,y)\ge a$ for $x,y\in G$ and for each $x\in M$, there exists $y\in G$ such that $d(x,y)\le b$.

\parindent=0mm \vspace{.1in}
{\bf Theorem 3.B.10 ([25]):} Let $X$ be a Banach space such that for each closed subspace $M$ of $X$, the quotient map $\varphi: X\to X/M$ admits a Lipschitz right inverse on a net $G$: there exists a Lipschitz map $\psi:G\to X$ such that  $\varphi\psi=\text{ identity map on}\, G$. Then $X$ is isomorphic to a Hilbert space.

\parindent=0mm \vspace{.1in}
As noted in Theorem 3.B.6, the existence of a global right inverse of the quotient map gives that $M$ is a Lipschitz retract. However, it can be shown that under the conditions of the theorem, $M$ is locally complemented in $X$ and the conclusion follows from the last assertion of Corollary 3.A.8. The full details of the proof shall appear elsewhere.

\parindent=0mm \vspace{.1in}

{\bf (C).	Grothendieck and Kwapien Theorems} 

\parindent=0mm \vspace{.1in}

{\bf Theorem 3.C.1([7]},(see also [23] for a different proof): For a Banach space $X$, TFAE:

\parindent=0mm \vspace{.1in}
(i). Both $X$ and $X^*$ are linearly embedded in an $L_1$ space.

\parindent=0mm \vspace{.1in}
(ii). $X$ is a Hilbert space.

\parindent=0mm \vspace{.1in}
{\bf Theorem 3.C.2 ([11]):}  Let $X$ be a Banach space such that $X$ is linearly isomorphically embeddable with a subspace of each of $L_q$ and $L_p$  where $1\le p\le 2\le q<\infty$. Then $X$ is (linearly and isomorphically) a Hilbert space.

\parindent=0mm \vspace{.1in}
The above theorem is a direct consequence of Kwapien’s theorem [12] to the effect that a Banach space having type 2 and cotype 2 is linearly isomorphic to a Hilbert space. 

\parindent=0mm \vspace{.1in}
We provide below the Lipschitz analogue of Theorems 3.C.1 and 3.C.2 which can be proved using Ribe’s theorem [22] combined with the machinery of ultra products as covered in [4], Chapter 8.

\parindent=0mm \vspace{.1in}

{\bf Theorem 3.C.3 ([25]):} For a Banach space $X$, TFAE:

\parindent=0mm \vspace{.1in}
(i). Both $X$ and $X^*$ are Lipschitz embedded in an $L_1$ space

\parindent=0mm \vspace{.1in}
(ii). $X$ is a Hilbert space.

\parindent=0mm \vspace{.1in}
Finally, a simple argument involving type and cotype of a Banach space gives the following Lipschitz analogue of Kwapien’s Theorem 3.C.2 which avoids the use of the highly non-trivial theorems of Mankiewicz [14] and Enflo [5]. 

\parindent=0mm \vspace{.1in}
{\bf Theorem 3.C.4 ([25]):} Let $X$ be a Banach space such that $X$ is Lipschitz embeddable with a subspace of each of $L_q$ and $L_p$ where $1\le p\le 2\le q<\infty$. Then $X$ is isomorphically a Hilbert space.

\parindent=0mm \vspace{.1in}
In contrast to Banach spaces as covered in the previous theorem, there are metric spaces $M$ (the infinite tree, for example) that are Lipschitz embeddable into the (cotype 2) Banach space $\ell_1$ and the James space $\mathbb J$- a nonre flexive Banach space of type 2 - but $M$ is not bi-Lipschitz embeddable in a Hilbert space. This was shown by A. Naor and R. Young in an important work [19] where they investigate, among other things, the question of extending Maurey’s factorisation theorem [17] and Theorem 3.C.2 stated above to the metric setting. In this connection, we quote directly from [19] on the impossibility of defining appropriate notions of type 2 and cotype 2 for metric spaces ``that are bi-Lipschitz invariant, pass to subsets, coincide for Banach spaces with type 2 and cotype 2, and for which Kwapien's theorem [12] holds for doubling metric spaces, i.e., any doubling space that has both type 2 and cotype 2 admits a bi-Lipschitz embedding into a Hilbert space.  The following theorem supplements these observations in the more familiar setting of $\ell_p$ spaces.

\parindent=0mm \vspace{.1in}
{\bf Theorem 3.C.5([19]):} For any $2< p<\infty$ there is a Lipschitz mapping $f:\ell_p\to \ell_1$  that cannot be factored through a subset of a Hilbert space using Lipschitz mappings. Also, there exists a metric space $M$ that Lipschitz embeds into both $\ell_1$ and $\ell_p$, but $M$ does not Lipschitz embed into a Hilbert space. 

\parindent=0mm \vspace{.1in}
{\bf (D). Norm attaining vectors and norm attaining functionals}

\parindent=0mm \vspace{.1in}
In recent years, the study of norming vectors and their dual counterpart involving norm attaining vectors in Banach spaces has received considerable attention, leading to new insights into the structure and geometry of Banach spaces. Of special interest has been the investigation of properties ensuring the existence of (high dimensional) vector spaces comprising exclusively of norming vectors and of norm attaining linear functionals. Here we quote just two results, one from each of these genres of norming objects as a sample from a vast repertoire of results which have been proved in recent years.

\parindent=0mm \vspace{.1in}
{\bf Theorem 3.D.1([1]):} The norm of a finite dimensional space $X$ arises from an inner product if and only if for each linear map $T:X\to X$, the set $N(T)$ of norming vectors of $T$ is a vector space, where $N(T)=\big\{x\in X:\|T(x)\|= \|T\|\|x\|\big\}$.

\parindent=0mm \vspace{.1in}
In what follows, for a given Banach space $X$, we shall use the symbol $NA(X)$ to denote the set of all norm-attaining functionals on $X: NA(X)=\big\{f\in X^*: x\in X \, \text{s.t.}\, \|f(x)\|=\|f\| \|x\|\big\}$.

\parindent=0mm \vspace{.1in}
As an easy consequence of the Hahn-Banach theorem, it follows that for an infinite dimensional dual Banach space $Z=X^*$, the set $NA(Z)$ contains ($X$ as) an infinite dimensional closed) subspace. However, there also exist non-dual Banach spaces enjoying this property; $NA(c_0)$ contains an infinite dimensional subspace as is also the space $L_1(\mu)$ where $\mu$ is a $\sigma$-finite measure, provided the latter space is infinite dimensional. On the other hand, there also exist infinite dimensional Banach spaces $X$ with $NA(X)$ not containing even a 2-dimensional subspace [21].

\parindent=0mm \vspace{.1in}
On the other hand, it's a famous theorem of R. C. James that for a Banach space $X$, the equality $X^*=NA(X)$ characterizes $X$ as a reflexive Banach space. In this situation, a natural question arises whether for each $f\in X^*$, it is possible to choose $x\in X$ in a linear isometric manner such that $\|f\|\|x\|=|f(x)|$. As a special case, this motivates the question of describing Banach spaces which are linearly isometric with their duals.  As an important theme in Banach space theory, we note that whereas Hilbert spaces provide easy examples of such spaces, converse doesn’t always hold. This is testified by taking X to be a non-Hilbertian reflexive Banach space and note that the map   $\varphi : x\times X^*\to \big(x\times X^*\big)^*$ given by $\varphi(x,y)=(y,x)$ defines a linear isometry. As a finite dimensional example, one may consider the space $X=\ell_1^2$ and its dual space $X^*=\ell_\infty^2$ 

\parindent=0mm \vspace{.1in}
The above discussion suggests that certain additional conditions have to be imposed on the linear isometry to ensure that $X$ is a Hilbert space. A remarkable result due to Mascioni [15] says that if for a Banach space $X$ there exists $c>0$ such that every infinite dimensional subspace of $\ell_2\{X\}$ is $c$-isomorphic with its dual, then $X$ is (isomorphic to) a Hilbert space.

\parindent=0mm \vspace{.1in}
On the other hand, one may impose extra conditions on the linear map $\varphi:X^*\to X$ to ensure that $X$ is (isomorphically) a Hilbert space. One  such condition involves the requirement that $\langle \varphi (x^*), x^*\rangle =\|x^*\|^2$, for all $x^*\in X^*$, in which case $X$ is (isometrically) a Hilbert space. More generally, the existence of $c>0$ such that $\langle \varphi (x^*), x^*\rangle \ge c\|x^*\|^2$, for all $x^*\in X^*$ yields that $X$ is linearly isomorphic to a Hilbert space. In this case, the inner product on $X^*$ may be defined by 
$$\langle x^*, y^*\rangle= \dfrac{1}{2} \Big\{\langle \varphi (x^*), y^*\rangle+ \langle \varphi (y^*), x^*\rangle\Big\}$$
which is easily seen to yield a norm equivalent to the dual norm on $X^*$. Consequently, $X^*$ and hence $X$ is isomorphic to a Hilbert space. 

\parindent=0mm \vspace{.1in}
Existence of Lipschitz selection of norm attaining functionals

\parindent=0mm \vspace{.1in}
In the next theorem, we address the question of the existence of a Lipschitz analogue of the foregoing statement (Section 3) and show how the existence of a selection map involving a bi-Lipschitz choice of norm-one vectors in a Banach space $X$ for the set of all bounded linear functionals on $X$ places severe restrictions on $X$. To the best of our understanding, the author has come across no reference in the literature pertaining to a characterisation of Hilbert spaces in terms of the Lipschitz structure of norm attaining functionals. More precisely, we have

\parindent=0mm \vspace{.1in}
{\bf Theorem 3.D.2 ([25]):} Let $X$ be a Banach space and $\varphi:X^*\to X$ a bilipschitz isomorphism such that there exists $c>0$ with
$$  \langle \varphi (x^*), x^*\rangle\ge c\|x^*\|, ~\forall ~x^*\in X^*.$$     
Then $X$ is isomorphic to a Hilbert space, and conversely.

\end{document}